\documentclass[reqno]{amsart}  
\usepackage[T1]{fontenc}
\usepackage{amsxtra}

\usepackage{color}
\usepackage{amscd,amsfonts,amsmath,amssymb,
amsthm,amstext,latexsym}
\usepackage{pifont,graphicx}
\usepackage{enumerate}
\usepackage[english]{babel}
\usepackage[all,cmtip]{xy} % for commutative diagram
\DeclareFontFamily{U}{mathx}{\hyphenchar\font45}
\DeclareFontShape{U}{mathx}{m}{n}{
      <5> <6> <7> <8> <9> <10>
      <10.95> <12> <14.4> <17.28> <20.74> <24.88>
      mathx10
      }{}
\DeclareSymbolFont{mathx}{U}{mathx}{m}{n}
\DeclareFontSubstitution{U}{mathx}{m}{n}
\DeclareMathAccent{\widecheck}{0}{mathx}{"71}
\DeclareMathAccent{\widetilde}{0}{mathx}{"72}
\DeclareMathAccent{\widebar}{0}{mathx}{"73}
\DeclareMathAccent{\widevec}{0}{mathx}{"74}
\DeclareMathAccent{\widehat}{0}{mathx}{"70}

\newcommand{\A}{{\mathbb A}}
\newcommand{\Z}{{\mathbb Z}}
\renewcommand{\S}{{\mathbb S}}
\newcommand{\R}{{\mathbb R}}
\newcommand{\C}{{\mathbb C}}
\newcommand{\D}{{\mathbb D}}

\renewcommand{\H}{{\mathbb H}}

\def\boC{{\mathcal C}}
\def\boD{{\mathcal D}}

\def\boF{{\mathcal F}}
\def\boG{{\mathcal G}}
\def\boH{{\mathcal H}}

\def\boM{{\mathcal M}}

\def\boW{{\mathcal W}}
\def\cqfd{\hfill$\Box$}
\def\Res{{\,\rm Res}}

\def\Re{{\rm Re}}
\def\Im{{\rm Im}}
\def\wcPhi{{\widecheck{\Phi}}}
\def\wcxi{{\widecheck{\xi}}}
\def\wcG{{\widecheck{G}}}

\def\wtPhi{{\widetilde{\Phi}}}

\def\wtf{{\widetilde{f}}}

\def\wtM{{\widetilde{M}}}
\def\wtN{\widetilde{N}}
\def\wtOmega{{\widetilde{\Omega}}}

\def\wtgamma{{\widetilde{\gamma}}}

\def\wtSigma{\widetilde{\Sigma}}
\def\wtz{\widetilde{z}}
\def\wtU{\widetilde{U}}
\def\wtV{\widetilde{V}}
\def\wtpsi{\widetilde{\psi}}
\def\whPhi{{\widehat{\Phi}}}

\def\whf{{\widehat{f}}}
\def\whxi{{\widehat{\xi}}}

\def\whomega{\widehat{\omega}}
\def\ii{{\rm i}}
\def\sl{\mathfrak{sl}}
\def\su{\mathfrak{su}}

\def\Uni{{\rm Uni}}
\def\Pos{{\rm Pos}}
\def\Sym{{\rm Sym}}
\def\Nor{{\rm Nor}}
\def\x{{\bf x}}
\def\y{{\bf y}}
\newcommand{\smallfrac}[2]{\mbox{$\frac{#1}{#2}$}}
\renewcommand{\matrix}[1]{\left(\begin{array}{cc} #1\end{array}\right)}
\newcommand{\minimatrix}[1]{\left(\begin{smallmatrix}#1\end{smallmatrix}\right)}
\newtheorem{theorem}{Theorem}

\newtheorem{proposition}{Proposition}
\newtheorem{remark}{Remark}
\newtheorem{corollary}{Corollary}
\newtheorem{claim}{Claim}

\newtheorem{definition}{Definition}

\title{Construction of constant mean curvature $n$-noids using the DPW method}
\author{Martin Traizet}
\begin{document}
\maketitle
{\em Abstract: we construct constant mean curvature surfaces in euclidean space with 
genus zero and $n$
ends asymptotic to Delaunay surfaces using the DPW method.}
\bigskip
\section{Introduction}
\label{introduction}
In \cite{dorfmeister-pedit-wu}, Dorfmeister, Pedit and Wu have shown that harmonic maps
from a Riemann surface to a symmetric space admit a Weierstrass-type
representation, which means that they can be represented in terms of
holomorphic data.
In particular, surfaces with constant mean curvature one (CMC-1 for short) in
euclidean space admit such a representation, owing to the fact that the Gauss map
of a CMC-1 surface is a harmonic map to the 2-sphere.
This representation is now called the DPW method and has been widely used
to construct CMC-1 surfaces in $\R^3$ and also constant mean curvature surfaces in homogeneous spaces such as the sphere $\S^3$ or hyperbolic space $\H^3$:
see for example \cite{dorfmeister-haak,
dorfmeister-wu,
heller-heller-schmitt,
heller1,
heller2,
heller3,
kilian-kobayashi-rossman-schmitt,
kilian-mcintosh-schmitt,
schmitt,
schmitt-kilian-kobayashi-rossman}.
Also the DPW method has been implemented by N. Schmitt to make computer
images of CMC-1 surfaces.
\medskip

The main limitation to the construction of examples is the Monodromy
Problem,
so either the topology of the constructed examples is limited
or symmetries are imposed to the construction, in order to reduce the number of equations to be solved.
\medskip

In constract, Kapouleas \cite{kapouleas} has constructed embedded CMC-1 surfaces with
no limitation on the genus or number of ends by gluing round spheres and
pieces of Delaunay surfaces, using Partial Differential Equations techniques.
An interesting question is wether similar results can be achieved with the DPW
method.
In this paper, we make a first step in this direction by constructing $n$-noids:
genus zero CMC-1 surfaces with $n$ ends.
\begin{theorem}
\label{theorem1}
Given $n\geq 3$ distinct unit vectors $u_1,\cdots,u_n$ in $\R^3$ and $n$ non-zero real weights $\tau_1,\cdots,
\tau_n$ satisfying the balancing condition 
$$\sum_{i=1}^n \tau_i u_i=0$$
there exists a smooth $1$-parameter family of CMC-1 surfaces $(M_t)_{0<t<\varepsilon}$
with genus zero, $n$ Delaunay ends and the following properties:
\begin{enumerate}[1.]
\item If we denote $w_{i,t}$ the weight of the $i$-th Delaunay end and $\Delta_{i,t}$ its
axis, then $$\lim_{t\to 0}\frac{w_{i,t}}{t}=8\pi\tau_i$$ and $\Delta_{i,t}$ converges to
the half-line through the origin directed by $u_i$.
\item If all weights $\tau_i$ are positive, then $M_t$ is Alexandrov-embedded.
\item If moreover
the angle between $u_i$ and $u_j$ is greater than $\frac{\pi}{3}$ for all $j\neq i$, then
$M_t$ is embedded.
\end{enumerate}
\end{theorem}
\begin{figure}
\begin{center}
\includegraphics[width=4cm]{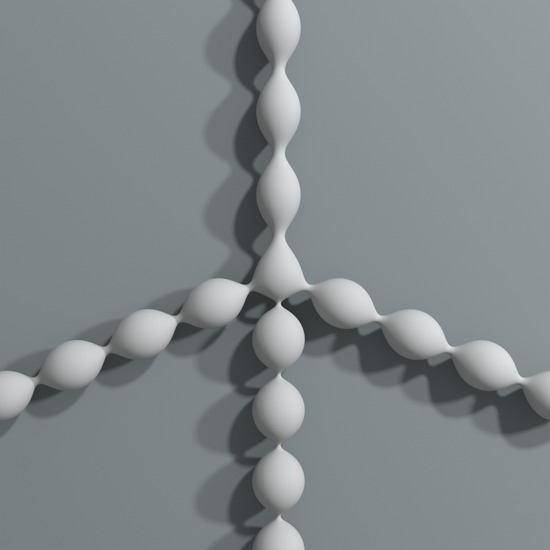}
\end{center}
\label{figure1}
\caption{a symmetric 4-noid. Image by N. Schmitt}
\end{figure}
These examples can be described heuristically as the unit sphere with $n$ half Delaunay
surfaces with small necksizes attached at the points $u_1,\cdots,u_n$
(see Figure \ref{figure1}).
They are a particular case of the construction of Kapouleas \cite{kapouleas}.
Here are some related results:
\begin{enumerate}[1.]
\item The full family of 3-noids is constructed in \cite{dorfmeister-wu,schmitt-kilian-kobayashi-rossman}
using the DPW method. Their construction is more general than ours
in the case $n=3$, since we only construct $3$-noids whose ends have small necksizes.
\item Highly symmetric $n$-noids with $n\geq 4$ have been constructed using the DPW method in \cite{kilian-mcintosh-schmitt,schmitt}.
\item DPW potentials of general $n$-noids are investigated in \cite{dorfmeister-wu,gerding-pedit-schmitt}, but the Monodromy Problem is not solved when $n\geq 4$. Quoting
\cite{gerding-pedit-schmitt}:
{\em Outside of experimental examples, the construction of $k$-noids for $k\geq 4$, even
under the additional assumption of coplanarity, remains elusive from the integrable
systems perspective.}
\item Alexandrov-embedded $3$-noids have been classified in \cite{karsten-kusner-sullivan},
and Alexandrov-embedded $n$-noids with coplanar ends have been classified
in \cite{karsten-kusner-sullivan2} (with methods unrelated to DPW).
\item The balancing condition of Theorem \ref{theorem1} is necessary by the general balancing
formula for CMC-1 surfaces (Theorem 3.4 in \cite{korevaar-kusner-solomon}).
\end{enumerate}
\medskip

Our goal in this paper is to prove Theorem \ref{theorem1} with the DPW method, using a rather simple and natural DPW potential,  inspired from the one used in \cite{schmitt-kilian-kobayashi-rossman} for $3$-noids.
Our main motivation is to make progress in the DPW method. In this regard, here
is what is achieved in this paper:
\begin{enumerate}[1.]
\item We are able to solve the Monodromy Problem on an $n$-punctured sphere,
whose fundamental group has $n-1$ generators, without any symmetry assumption.
The Monodromy Problem is solved by an Implicit Function argument in a suitable
functional space.
\item We are able to prove embeddedness. I believe this is the first time that a non-trivial example is proven to be embedded using the DPW method.
The proof relies on the study of Delaunay ends in \cite{kilian-rossman-schmitt}
and \cite{raujouan}.
\item The DPW potential that we use has the property that its poles depend on
the spectral parameter $\lambda$. This  problem is answered by Theorem \ref{theorem-change-coordinate} in
Section \ref{lambda-dependent}, a general result which allows $\lambda$-dependent changes of variable
in the DPW method. Theorem \ref{theorem-change-coordinate} adds some flexibility to the choice of the DPW potential.
\end{enumerate}
\section{Background}
\label{background}
In this section, we recall standard notations and results used in the
DPW method.
We work in the ``untwisted'' setting.
For a comprehensive introduction to the DPW method, we suggest
\cite{fujimori-kobayashi-rossman}.
\subsection{Loop groups}
\label{background-section1}
A loop is a smooth map from the unit circle $\S^1=\{\lambda\in\C:|\lambda|=1\}$ to a matrix group.
The circle variable is denoted $\lambda$ and called the spectral parameter.
For $\rho>0$, we denote $\D_{\rho}=\{\lambda\in\C\,:\, |\lambda|<\rho\}$,
$\D_{\rho}^*=\D_{\rho}\setminus\{0\}$ and $\D=\D_1$.
\begin{itemize}
\item If $G$ is a matrix Lie group (or Lie algebra), $\Lambda G$ denotes the group
(or algebra)
of smooth maps $\Phi:\S^1\to G$.
%\item $\Lambda SL(2,\C)$ is the set of smooth maps $\Phi:\S^1\to SL(2,\C)$.
%\item $\Lambda SU(2)\subset\Lambda SL(2,C)$ is the set of smooth maps $F:\S^1\to SU(2)$.
\item $\Lambda_+ SL(2,\C)\subset \Lambda SL(2,\C)$ is the subgroup of smooth maps $B:\S^1\to SL(2,\C)$ which extend holomorphically to $\D$.
\item $\Lambda_+^{\R}SL(2,\C)$ is the subgroup of $B\in \Lambda_+ SL(2,\C)$ such that
$B(0)$ is upper triangular with positive elements on the diagonal.
\end{itemize}
\begin{theorem}[Iwasawa decomposition]
The multiplication $\Lambda SU(2)\times \Lambda_+^{\R} SL(2,\C)\to\Lambda SL(2,\C)$
is a diffeomorphism. The unique splitting of an element $\Phi\in\Lambda SL(2,\C)$ as
$\Phi=FB$ with $F\in\Lambda SU(2)$ and $B\in\Lambda_+^{\R} SL(2,\C)$
is called Iwasawa decomposition.
$F$ is called the unitary factor of $\Phi$ and
denoted $\Uni(\Phi)$.
$B$ is called the positive factor and denoted $\Pos(\Phi)$.
\end{theorem}
\subsection{The matrix model of $\R^3$}
\label{background-section2}
In the DPW method, one identifies $\R^3$ 
with the Lie algebra $\su(2)$ by
$$x=(x_1,x_2,x_3)\in\R^3\longleftrightarrow X=\frac{-\ii}{2}\matrix{
-x_3&x_1+\ii x_2\\x_1-\ii x_2 &x_3}\in\mathfrak{su}(2).$$
(This is essentially the same as identifying $\R^3$ with imaginary
quaternions and using a matrix model for quaternions.)
Under this identification, the euclidean scalar product and norm are given by
$$\langle x,y\rangle=-2\mbox{tr}(XY),\qquad
||x||^2=4\det(X).$$
The group $SU(2)$ acts as linear isometries on $\su(2)$ by
$H\cdot X=HXH^{-1}$.
%The kernel of this action is $\pm I_2$ so that
%$SO(3)$ is isomorphic to $SU(2)/\{\pm I_2\}$.
\subsection{The DPW method}
\label{background-section3}
The input data for the DPW method is a quadruple $(\Sigma,\xi,z_0,\phi_0)$ where:
\begin{itemize}
\item $\Sigma$ is a Riemann surface.
\item
$\xi=\xi(z,\lambda)$ is a $\Lambda\sl(2,\C)$-valued holomorphic 1-form on $\Sigma$
called the DPW potential. More precisely,
\begin{equation}
\label{equation-xi}
\xi=\matrix{\alpha&\lambda^{-1}\beta\\ \gamma&-\alpha}
\end{equation}
where $\alpha(z,\lambda)$, $\beta(z,\lambda)$, $\gamma(z,\lambda)$
are holomorphic 1-forms on $\Sigma$ with respect to the $z$ variable,
and are holomorphic with respect to $\lambda$ in the disk $\D_{\rho}$
for some $\rho>1$.
\item $z_0\in\Sigma$ is a base point.
\item $\phi_0\in\Lambda SL(2,\C)$ is an initial condition.
\end{itemize}
Given this data, the DPW method is the following procedure.
Let $\wtSigma$ be the universal cover of $\Sigma$ and $\wtz_0\in\wtSigma$
be an arbitrary element in the fiber of $z_0$.
\begin{enumerate}[1.]
\item Solve the Cauchy Problem on $\wtSigma$:
$$d_z\Phi(z,\lambda)=\Phi(z,\lambda)\xi(z,\lambda)$$
 with initial condition
$$\Phi(\wtz_0,\lambda)=\phi_0(\lambda)$$
to obtain a solution
$\Phi: \wtSigma\to\Lambda SL(2,\C)$.
(The notation $d_z$ means that we are considering the differential with respect
to the $z$-variable. The lift of $\xi$ to $\wtSigma$ is still denoted $\xi$.)
\item Compute, for $z\in\wtSigma$, the unitary part $F(z,\cdot)=\Uni(\Phi(z,\cdot))$ in the Iwasawa decomposition of $\Phi(z,\cdot)$.
It is known that $F(z,\lambda)$ depends real-analytically on $z$.
\item Define $f:\wtSigma\to\su(2)\sim\R^3$ by the Sym-Bobenko formula:
$$f(z)=\ii\frac{\partial F}{\partial\lambda}(z,1)F(z,1)^{-1}=:\Sym(F(z,\cdot)).$$
Then $f$ is a CMC-1 (branched) conformal immersion. Its Gauss map is given by
$$N(z)=\frac{-\ii}{2}F(z,1)\matrix{1 &0 \\ 0 & -1}F(z,1)^{-1}=:\Nor(F(z,\cdot)).$$
Moreover, $f$ is regular at $z$ (meaning unbranched) if and only if $\beta(z,0)\neq 0$.
\end{enumerate}
Observe that we require the potential $\xi$ to be holomorphic on $\Sigma$.
In many examples, $\Sigma$ is a compact Riemann surface $\overline{\Sigma}$ minus a finite number of points, and $\xi$ extends meromorphically at the punctures.
\subsection{The Monodromy Problem}
Assume that $\Sigma$ is not simply connected so its universal cover $\wtSigma$
is not trivial.
Let $\mbox{Deck}(\wtSigma/\Sigma)$ be the group of fiber-preserving diffeomorphisms of $\wtSigma$.
For $\gamma\in\mbox{Deck}(\wtSigma/\Sigma)$, let
$$\boM_{\gamma}(\Phi)(\lambda)=\Phi(\gamma(z),\lambda)\Phi(z,\lambda)^{-1}$$ 
be the monodromy of $\Phi$ with respect to $\gamma$
(which is independent of $z\in\wtSigma$).
The standard condition which ensures that the immersion $f$ descends to
a well defined immersion on $\Sigma$
is the following system of equations, called the Monodromy Problem.
\begin{equation}
\label{monodromy-problem}
\forall \gamma\in\mbox{Deck}(\wtSigma/\Sigma)\quad
\left\{\begin{array}{lc}
\boM_{\gamma}(\Phi)\in\Lambda SU(2)\qquad &(i)\\
\boM_{\gamma}(\Phi)(1)=\pm I_2\qquad&(ii)\\
\frac{\partial\boM_{\gamma}(\Phi)}{\partial\lambda}(1)=0\qquad &(iii)
\end{array}\right.\end{equation}
Indeed, Condition (i) implies that $F$ has a monodromy, equal to the monodromy
of $\Phi$, and provided this is true,
Conditions (ii) and (iii) are equivalent to the fact that $f$ is well defined on $\Sigma$.
\medskip

One can identify $\mbox{Deck}(\wtSigma/\Sigma)$ with the fundamental group
$\pi_1(\Sigma,z_0)$ (see for example Theorem 5.6 in \cite{forster}), so we will in general see $\gamma$ as an element of
$\pi_1(\Sigma,z_0)$. This identification, however, is not canonical, as it depends on
the choice of $\wtz_0$. Under this identification, the monodromy of $\Phi$ with
respect to $\gamma\in\pi_1(\Sigma,z_0)$ is given by
$$\boM_{\gamma}(\Phi)(\lambda)=\Phi(\wtgamma(1),\lambda)\Phi(\wtgamma(0),\lambda)^{-1}$$
where $\wtgamma:[0,1]\to\wtSigma$ is the lift of $\gamma$ such that $\wtgamma(0)=\wtz_0$.
\subsection{Basic examples}
\label{section-basic-examples}
\begin{enumerate}[1.]
\item A round sphere is obtained with the data
$$\Sigma=\C,
\qquad
\xi(z,\lambda)=\matrix{0&\lambda^{-1}\\0&0}dz,
\qquad
z_0=0,
\qquad
\phi_0=I_2.$$
The solution of the Cauchy Problem is
$$\Phi(z,\lambda)=\matrix{1&\lambda^{-1}z\\0&1}.$$
Its Iwasawa decomposition is
\begin{equation}
\label{Iwasawa-sphere}
F(z,\lambda)=\frac{1}{\sqrt{1+|z|^2}}\matrix{1&\lambda^{-1}z\\-\lambda\overline{z}&1},
\qquad
B(z,\lambda)=\frac{1}{\sqrt{1+|z|^2}}\matrix{1&0\\\lambda\overline{z}&1+|z|^2}.
\end{equation}
The Sym-Bobenko formula gives
$$f(z)=\frac{1}{1+|z|^2}\left(2\,\Re(z),2\,\Im(z),-2|z|^2\right)=(0,0,-1)+\pi^{-1}(z)$$
where
\begin{equation}
\label{inverse-stereographic}
\pi^{-1}(z)=\left(\frac{2\,\Re(z)}{1+|z|^2},\frac{2\,\Im(z)}{1+|z|^2},\frac{1-|z|^2}{1+|z|^2}\right)
\end{equation}
is the inverse stereographic projection from the south pole.
The immersion $f$ extends smoothly at $\infty$ and gives a conformal parametrisation
of the unit sphere centered at $(0,0,-1)$.
\item Delaunay surfaces are obtained with the data
$$\Sigma=\C\setminus\{0\},
\qquad
\xi(z,\lambda)=\matrix{0&r\lambda^{-1}+s\\r\lambda+s&0}\frac{dz}{z},
\qquad
z_0=1,
\qquad
\phi_0=I_2$$
where $r,s$ are non-zero real numbers such that $r+s=\frac{1}{2}$.
\end{enumerate}
\subsection{Gauging}
\label{background-section6}
\begin{definition}
A gauge on $\Sigma$ is a map $G:\Sigma\to\Lambda_+ SL(2,\C)$ such that
$G(z,\lambda)$ depends holomorphically on $z\in\Sigma$ and
$\lambda\in\D_{\rho}$ and $G(z,0)$
is upper triangular (with no restriction on its diagonal elements).
\end{definition}
Let $\Phi$ be a solution of $d_z\Phi=\Phi\xi$ and $G$ be a gauge.
Let $\whPhi=\Phi\times G$. Then $\whPhi$ and $\Phi$
define the same immersion $f$. This is called ``gauging''.
The gauged potential is
$$\whxi=\whPhi^{-1}d_z\whPhi=G^{-1}\xi G+G^{-1}d_z G$$
and will be denoted $\xi\cdot G$, the dot denoting the action of the gauge
group on the potential.
\section{Change of variable depending on $\lambda$}
\label{lambda-dependent}
The following theorem allows us to use $\lambda$-dependent changes of variable
in the DPW method.
We will use this theorem in Section \ref{section-delaunay-ends}
with $\psi_{\lambda}$ a translation depending on $\lambda$.
\medskip

Let $U$ and $V$ be Riemann surfaces, $\rho>1$ and
$\psi:U\times\D_{\rho}\to V$ be a holomorphic map.
For $(z,\lambda)\in U\times\D_{\rho}$, we write $\psi_{\lambda}(z)=\psi(z,\lambda)$.
Let $\xi(z,\lambda)$ be a DPW potential on $V$.
\begin{theorem}
\label{theorem-change-coordinate}
Assume that $V$ is simply connected.
Let $\Phi(z,\lambda)$ be a solution of $d_z\Phi=\Phi\xi$, holomorphic on
$V\times\D_{\rho}^*$,
and $f:V\to\R^3$ be the immersion defined by $f=\Sym(\Uni(\Phi))$.
Define:
$$\whPhi(z,\lambda)=\Phi(\psi_{\lambda}(z),\lambda)\quad\mbox{ for } (z,\lambda)\in U\times\D_{\rho}^*.$$
Let $\whf: U\to\R^3$ be the (branched) immersion defined by $\whf=\Sym(\Uni(\whPhi))$.
Then $\whf=f\circ\psi_0$ in $U$.
\end{theorem}
\begin{remark}
\label{remark-change-coordinate}

\begin{enumerate}[1.]
\item We relax the hypothesis that $V$ is simply connected in Corollary \ref{corollary-change-coordinate}.
\item $\whPhi$ solves $d\whPhi=\whPhi\whxi$ in $U$ where $\whxi(\cdot,\lambda)
=\psi_{\lambda}^*\xi(\cdot,\lambda)$.
\item
Let $(F,B)$ be the Iwasawa decomposition
of $\Phi$.
Then of course
$$\whPhi(z,\lambda)=F(\psi_{\lambda}(z),\lambda)\times B(\psi_{\lambda}(z),\lambda)$$
but this is in general not the Iwasawa decomposition of $\whPhi(z,\lambda)$.
It is true that $F(\psi_{\lambda}(z),\lambda)\in\Lambda SU(2)$.
But in general, $B(z,\lambda)$ does not depend holomorphically on $z$, so
there is no reason that $B(\psi_{\lambda}(z),\lambda)$ should extend holomorphically
to $\lambda\in\D$.
For example, in the case of the spherical potential (see Section \ref{section-basic-examples}),
if $B(\psi_{\lambda}(z),\lambda)$ is holomorphic with respect to $\lambda$,
Equation \eqref{Iwasawa-sphere} gives that
$|\psi_{\lambda}(z)|$ is holomorphic, hence $\psi_{\lambda}(z)$
is constant with respect to $\lambda$.
\end{enumerate}
\end{remark}
Proof of Theorem \ref{theorem-change-coordinate}: define for $z\in U$ and $\lambda\in \D_{\rho}^*$:
$$\wtPhi(z,\lambda)=\Phi(\psi_0(z),\lambda).$$
Then
$$\wtPhi(z,\lambda)=F(\psi_0(z),\lambda)\times B(\psi_0(z),\lambda)$$
is the Iwasawa decomposition of $\wtPhi(z,\lambda)$.
(The difference with Point 3 of Remark \ref{remark-change-coordinate} is that $\psi_0(z)$ does not depend on $\lambda$.)
Let $\wtf=\Sym(\Uni(\wtPhi))$.
By the Sym-Bobenko formula, we obtain for $z\in U$
$$\wtf(z)=f(\psi_0(z)).$$
Define for $z\in U$ and $\lambda\in\D_{\rho}^*$:
$$G(z,\lambda)=\wtPhi(z,\lambda)^{-1}\times\whPhi(z,\lambda)
=\Phi(\psi_0(z),\lambda)^{-1}\times\Phi(\psi_{\lambda}(z),\lambda).$$
The following two claims prove that $G$ is a gauge.
Hence $\wtPhi$ and $\whPhi$
are gauge-equivalent, so $\wtf=\whf$ in $U$. This proves Theorem \ref{theorem-change-coordinate}.\cqfd
\begin{claim}
\label{claim2}
$G(z,\lambda)$ extends holomorphically at $\lambda=0$.
\end{claim}
Proof: we use Gr\"onwall inequality to estimate $G(z,\lambda)$.
Fix $z\in U$ and $\lambda\in \D^*$. Consider the path $\gamma:[0,1]\to V$ defined by $\gamma(s)=\psi_{s\lambda}(z)$.
Then for $s\in[0,1]$
$$\frac{d}{ds}\Phi(\gamma(s),\lambda)=d_z\Phi(\gamma(s),\lambda)\,\gamma'(s)=\Phi(\gamma(s),\lambda)\,
\xi(\gamma(s),\lambda)\,\gamma'(s).$$
Hence for $t\in[0,1]$
$$\Phi(\gamma(t),\lambda)=\Phi(\gamma(0),\lambda)+\int_0^t
\Phi(\gamma(s),\lambda)\,
\xi(\gamma(s),\lambda)\,\gamma'(s)\,ds.$$
Multiplying by $\Phi(\gamma(0),\lambda)^{-1}$ on the left an taking norms, we obtain
$$||\Phi(\gamma(0),\lambda)^{-1}\,\Phi(\gamma(t),\lambda)||
\leq 1+\int_{0}^t||\Phi(\gamma(0),\lambda)^{-1}\,\Phi(\gamma(s),\lambda)||\cdot
||\xi(\gamma(s),\lambda)\gamma'(s)||\,ds.$$
By Gr\"onwall inequality, we obtain
\begin{equation}
||G(z,\lambda)||=||\Phi(\gamma(0),\lambda)^{-1}\,\Phi(\gamma(1),\lambda)||
\leq\exp\int_0^1||\xi(\gamma(s),\lambda)\gamma'(s)||\,ds.
\end{equation}
We have
$$\gamma'(s)=\lambda\frac{\partial\psi}{\partial\lambda}(z,s\lambda).$$
Hence
$$||\xi(\gamma(s),\lambda)\gamma'(s)||
=\left\|\lambda\xi(\psi_{s\lambda}(z),\lambda)
\frac{\partial\psi}{\partial\lambda}(z,s\lambda)\right\|\leq c$$
for some constant $c$ independent of $\lambda\in\D^*$, because $\xi$ has
a simple pole at $\lambda=0$ and using continuity. 
Hence $G(z,\lambda)$ is bounded for $\lambda\in\D^*$.
By Riemann extension theorem, $G(z,\lambda)$ extends holomorphically at $\lambda=0$.
\cqfd
\begin{claim}
\label{claim3}
$G(z,0)$ is upper triangular.
\end{claim}
Proof:
Define
$$K(\lambda)=\matrix{\frac{1}{\sqrt{\lambda}}&0\\0&\sqrt{\lambda}}$$
and let $\wcPhi=\Phi K$.
Then $\wcPhi$ solves $d_z\wcPhi=\wcPhi\wcxi$ where
$$\wcxi=\xi\cdot K=K^{-1}\xi K.$$
Of course, $K$ is not an admissible gauge and $\wcxi$ is not a DPW potential
(which does not matter here). But if
we write $\xi$ as in Equation \eqref{equation-xi},
we have
$$\wcxi=\matrix{\alpha&\beta\\\lambda^{-1}\gamma&-\alpha}$$
so $\wcxi$ has (at most) a simple pole at $\lambda=0$.
Define
$$\wcG(z,\lambda)=\wcPhi(\psi_0(z),\lambda)^{-1}\times\wcPhi(\psi_\lambda(z),\lambda).$$
By the proof of Claim \ref{claim2}, since $\wcxi$ has a simple pole at $\lambda=0$,
 $\wcG(z,\lambda)$ extends holomorphically at $\lambda=0$.
Now $G$ and $\wcG$ are related by
$$G(z,\lambda)=K(\lambda)\wcG(z,\lambda)K(\lambda)^{-1}.$$
This gives $G_{21}(z,\lambda)=\lambda \wcG_{21}(z,\lambda)$, so
$G_{21}(z,0)=0$.
\cqfd
\medskip

Let us illustrate the proof of Theorem \ref{theorem-change-coordinate} in a case where one can compute explicitely
the gauge $G$.
Consider a DPW potential of the form
$$\xi(z,\lambda)=A(\lambda)\omega(z).$$
Then any solution $\Phi$ has the following form:
$$\Phi(z,\lambda)=\Phi(z_0,\lambda)\exp\left[A(\lambda)\int_{z_0}^z \omega\right].$$
The gauge $G$ is given by
$$G(z,\lambda)=\exp\left[A(\lambda)\int_{\psi_{0}(z)}^{\psi_{\lambda}(z)}\omega\right]$$
and it is straightforward to check that $G$ extends holomorphically at $\lambda=0$
with upper triangular value.
\medskip

We now relax the hypothesis that $V$ is simply connected.
Let $p:\wtU\to U$ and $q:\wtV\to V$ be the universal covers of $U$ and $V$.
Since $\wtU\times\D_{\rho}$ is simply connected, the holomorphic map
$\psi:U\times\D_{\rho}\to V$ lifts to a holomorphic map
$\wtpsi:\wtU\times\D_{\rho}\to\wtV$.
\begin{corollary}
\label{corollary-change-coordinate}
Let $\Phi(z,\lambda)$ be a solution of $d_z\Phi=\Phi\xi$ on $\wtV\times\D_{\rho}^*$.
Let $\wtf:\wtV\to\R^3$ be the immersion defined by $\wtf=\Sym(\Uni(\Phi))$.
Assume that $\Phi$ solves the Monodromy Problem, so $\wtf$ descends to
$f:V\to\R^3$.
Define
$$\whPhi(z,\lambda)=\Phi(\wtpsi_{\lambda}(z),\lambda)
\quad \mbox{ for }(z,\lambda)\in \wtU\times\D_{\rho}^*.$$
Let $\whf: \wtU\to\R^3$ be the (branched) immersion defined by $\whf=\Sym(\Uni(\whPhi))$.
Then $\whf=f\circ\psi_0\circ p$ in $\wtU$.
In other words, $\whf$ descends to $f\circ\psi_0$.
\end{corollary}
Proof: by Theorem \ref{theorem-change-coordinate}, we have $\whf=\wtf\circ\wtpsi_0$ in $\wtU$.
The conclusion follows from the following commutative diagram:
$$\xymatrix{
\wtU  \ar[d]^{p}\ar[r]^{\wtpsi_0}& \wtV \ar[d]^{q} \ar[dr]^{\wtf}\\
U \ar[r]^{\psi_0} & V\ar[r]^{f} & \R^3\\
}
$$
\cqfd
%\begin{remark} One can check that the gauge $G$ in the proof of Theorem \ref{theorem-change-coordinate}, which a priori is defined on $\wtU$,
%actually descends to a well defined gauge on $U$.
%\end{remark}
\section{Functional spaces}
\label{functional-spaces}
In the next section, we propose a DPW potential for $n$-noids.
The parameters in the definition of this potential are functions of $\lambda$.
Since we plan to use the Implicit Function Theorem, 
we need to introduce suitable functional spaces.
We decompose a function $f:\S^1\to\C$ in Fourier
series
$$f(\lambda)=\sum_{i\in\Z}f_i\lambda^i$$
Fix some $\rho>1$ and define
$$||f||=\sum_{i\in\Z}|f_i|\rho^{|i|}$$
Let $\boW$ be the space of functions $f$ with finite norm.
This is a Banach algebra (classically called the Wiener algebra when $\rho=1$).
Functions in $\boW$ extend holomorphically to the annulus $\frac{1}{\rho}<|\lambda|<\rho$.
\medskip
\def\pos{{\geq 0}}
\def\pospos{{+}}
\def\neg{{\leq 0}}
\def\negneg{{-}}

We define $\boW^{\pos}$, $\boW^{\pospos}$, $\boW^{\neg}$ and $\boW^{\negneg}$ as the subspaces of functions $f$ such that $f_i=0$
for $i<0$, $i\leq 0$, $i>0$ and $i\geq 0$, respectively.
Functions in $\boW^{\pos}$ extend holomorphically to the disk $\D_{\rho}$ and
satisfy $|f(\lambda)|\leq ||f||$ for all $\lambda\in\D_{\rho}$.
We write $\boW^0\sim\C$ for the subspace of constant functions, so we have
a direct sum $\boW=\boW^{\negneg}\oplus\boW^0\oplus\boW^{\pospos}$. A function $f$
will be decomposed as $f=f^{\negneg}+f^0+f^{\pospos}$ with
$(f^-,f^0,f^+)\in\boW^{\negneg}\times\boW^0\times\boW^{\pospos}$.
\medskip

We define the star operator by
$$f^*(\lambda)=\overline{f\left(\frac{1}{\overline{\lambda}}\right)}=\sum_{i\in\Z}\overline{f_{-i}}\lambda^i$$
The involution $f\mapsto f^*$ exchanges $\boW^{\pos}$ and $\boW^{\neg}$.
We have $\lambda^*=\lambda^{-1}$ and $c^*=\overline{c}$ if $c$ is a constant.
A function $f$ is real on the unit circle if and only if $f=f^*$.
\section{The DPW potential}
\label{nnoids-section1}
We now start the proof of Theorem \ref{theorem1}.
Without loss of generality, we assume (by a rotation) that all vectors $u_i$ are
non vertical. Let $\pi:\S^2\to\C\cup\{\infty\}$ be the stereographic projection from the south pole. For $i\in[1,n]$, we define $\pi_i=\pi(u_i)\in\C^*$ and we introduce three $\lambda$-dependent parameters $a_i$, $b_i$ and $p_i$ in the space $\boW^\pos$.
The collection of these parameters is denoted $\x=(a_i,b_i,p_i)_{1\leq i\leq n}\in(\boW^{\pos})^{3n}$. The parameter $\x$ is in a neighborhood of a (constant) central value which
we denote $\x_0$: the central value of $p_i$ is $\pi_i$, the central value of $a_i$ is
$\tau_i$, and we will compute
the central value of $b_i$ in Section \ref{nnoids-section5}.
We define a meromorphic 1-form $\omega_{\x}$ on
the Riemann sphere $\C\cup\{\infty\}$, depending on $\lambda$ and the parameter $\x$,
by
\begin{equation}
\label{eq-omega}
\omega_{\x}(z,\lambda)=\sum_{i=1}^n\left(\frac{a_i(\lambda)}{(z-p_i(\lambda))^2}+\frac{b_i(\lambda)}{z-p_i(\lambda)}\right)dz.
\end{equation}
For $t$ in a neighborhood of $0$ in $\R$, we define the meromorphic DPW potential $\xi_{t,\x}$ by
\begin{equation}
\label{eq-potential}
\xi_{t,\x}(z,\lambda)=\matrix{0 & \lambda^{-1}dz\\ t(\lambda-1)^2\omega_{\x}(z,\lambda) & 0}
\end{equation}
and take the initial condition $z_0=0$, $\phi_0=I_2$.
Since the parameters $a_i$, $b_i$ and $p_i$ are holomorphic functions of $\lambda$
in the disk $\D_{\rho}$, $\xi_{t,\x}$ is an admissible DPW potential.
Here are some of its properties:
\begin{enumerate}[1.]
\item If $t=0$, we get the standard DPW data for the sphere (see Section \ref{section-basic-examples}). Therefore, for small $t\neq 0$,
we are constructing, away from the poles, a perturbation of the unit sphere centered at $(0,0,-1)$.
We could of course have chosen the initial condition $\phi_0$ so that the sphere is
centered at the origin: it suffices to take
$$\phi_0(\lambda)=\matrix{e^{\frac{\lambda^2-1}{4\lambda}}&0\\0&e^{\frac{1-\lambda^2}{4\lambda}}}\in\Lambda SU(2)$$
but it is simpler to take $\phi_0=I_2$ and translate afterwards the immersion by
$(0,0,1)$.
\item Thanks to the factor $(\lambda-1)^2$ in front of $\omega_{\x}$ in the definition of $\xi_{t,\x}$,
Equations (ii)  and (iii) of the Monodromy Problem \eqref{monodromy-problem} 
are automatically solved.
\item We will see in Section \ref{section-delaunay-ends} that in a neighborhood of $p_i$, the potential $\xi_{t,\x}$ is gauge-equivalent to a potential with a simple pole, yielding a fuchsian system of differential equations. Moreover, provided the Monodromy Problem is solved, the residue of the
gauged potential is a standard Delaunay residue.
Therefore, the immersion will have Delaunay ends by the work of
Kilian, Rossman and Schmitt in \cite{kilian-rossman-schmitt}.
\end{enumerate}
\begin{remark}
A potential of the form \eqref{eq-potential} has been introduced for 3-noids in \cite{schmitt-kilian-kobayashi-rossman}. The main difference is that in that paper, $\omega$ is fixed:
the parameter used to solve the Monodromy Problem is the
initial condition $\phi_0\in\Lambda SL(2,\C)$, but this works only for $n=3$.
\end{remark}
\section{The equations}
\subsection{Regularity at $\infty$}
\label{nnoids-section2}
The potential $\xi_{t,\x}$ has a double pole at $\infty$.
We want our immersion to extend smoothly at $\infty$, so we require the
potential to be gauge-equivalent to a regular potential at $\infty$.
Consider the gauge
$$G_{\infty}(z,\lambda)=\matrix{z & 0 \\ -\lambda & z^{-1}}.$$
The gauged potential is
$$\xi_{t,\x}\cdot G_{\infty}=\matrix{0 & \lambda^{-1}\frac{dz}{z^2}\\t(\lambda-1)^2z^2\omega_{\x} & 0}.$$
It is regular at $\infty$ if and only if $\omega_{\x}$ has a double zero at $\infty$.
Using the coordinate $w=1/z$ in a neighborhood of $\infty$, we obtain
$$\omega_{\x}=-\sum_{i=1}^n \left(b_iw^{-1}+(a_i+b_ip_i)w^0+(2a_ip_i+b_ip_i^2)w+O(w^2)\right)dw.$$
We define the following functions:
$$\boH_1(\x)=\sum_{i=1}^n b_i$$
$$\boH_2(\x)=\sum_{i=1}^n a_i+b_ip_i$$
$$\boH_3(\x)=\sum_{i=1}^n 2a_ip_i+b_ip_i^2.$$
We need to solve the equations $\boH_1(\x)=\boH_2(\x)=\boH_3(\x)=0$ so that $\omega_{\x}$
has a double zero at $\infty$.
\subsection{The Monodromy Problem}
\label{nnoids-section3}
Our potential has poles at $p_1,\cdots,p_n$ which are functions of $\lambda$.
Because the DPW method requires a fixed Riemann surface, we introduce
the following domain:
\begin{equation}
\label{def-Omega}
\Omega=\{z\in\C:\forall i\in[1,n],|z-\pi_i|>\varepsilon\}
\end{equation}
where $\varepsilon>0$ is a fixed, small enough number such that the disks
$D(\pi_i,8\varepsilon)$ for $1\leq i\leq n$ are disjoint and do not contain $0$.
We use the following standard notations for domains in the $z$-plane:
$$D(p,r)=\{z\in\C:|z-p|<r\}\quad\mbox{ and }\quad
D^*(p,r)=D(p,r)\setminus\{p\}.$$
For $\x$ close enough to $\x_0$, $\xi_{t,\x}$ is holomorphic
in $\Omega\times\D_{\rho}^*$.
Our goal is first to construct a family of immersion $f_t$ on $\Omega$.
Then we extend $f_t$ to an $n$-punctured Riemann sphere in Section \ref{section-delaunay-ends},
using Corollary \ref{corollary-change-coordinate}.
\medskip

Let $\wtOmega$ be the universal cover of $\Omega$
and $\Phi_{t,\x}(z,\lambda)$ be the solution of the Cauchy Problem
$d_z\Phi_{t,\x}=\Phi_{t,\x}\xi_{t,\x}$ on $\wtOmega$
with initial condition $\Phi_{t,\x}(\widetilde{0},\lambda)=I_2$.
For $i\in[1,n]$, we denote $\gamma_1,\cdots,\gamma_n$ a set of generators of the fundamental group $\pi_1(\Omega,0)$, with $\gamma_i$ encircling the point $\pi_i$
(in other words, freely homotopic in $\Omega$ to the circle $C(\pi_i,2\varepsilon)$).
Let
$$M_i(t,\x)=\boM_{\gamma_i}(\Phi_{t,\x})$$
be the monodromy of $\Phi_{t,\x}$ with respect to $\gamma_i$.
%By the theorem on smooth dependence on parameters for the solution of a linear system
%of differential equations, each element of $M_i$ is a smooth map from a neighborhood
%of $(0,\x_0)$ in
%$\R\times(\boW^{\pos})^3$ to $\boW$.
Provided the Regularity Problem at $\infty$ is solved, $\Phi_{t,\x}$ has no monodromy around $\infty$, so we need to solve the following Monodromy Problem
\begin{equation}
\label{monodromy-1}
M_i(t,\x)\in\Lambda SU(2)\quad
\mbox{ for $1\leq i\leq n-1$.}
\end{equation}
At $t=0$, we have $M_i(0,\x)=I_2$.
Recall that the matrix exponential is a local diffeomorphism from a neighborhood of $0$
in the Lie algebra $\sl(2,\C)$ (respectively $\su(2)$) to
a neighborhood of $I_2$ in $SL(2,\C)$ (respectively $SU(2)$).
The inverse diffeomorphism is denoted $\log$.
For $t\neq 0$ small enough and $\lambda\in\D_{\rho}\setminus\{1\}$, we define:
$$\wtM_i(t,\x)(\lambda):=\frac{\lambda}{t(\lambda-1)^2}\log M_i(t,\x)(\lambda).$$
Observe that $\lambda\in\S^1\Rightarrow\frac{(\lambda-1)^2}{\lambda}\in\R$.
So for $t\neq 0$ and $\lambda\neq 1$, Problem \eqref{monodromy-1}
is equivalent to the following Rescaled Monodromy Problem:
\begin{equation}
\label{monodromy-2}
\wtM_i(t,\x)\in\Lambda\su(2)
\quad\mbox{ for $1\leq i\leq n-1$.}
\end{equation}
The elements of the matrix $\wtM_i$ are denoted $\wtM_{i;k\ell}$ for $1\leq k,\ell\leq 2$.
\begin{proposition}
\label{proposition-Mi}
$\wtM_i(t,\x)(\lambda)$ extends smoothly at $t=0$ and $\lambda=1$, and
for $1\leq k,\ell\leq 2$, $\wtM_{i;k\ell}$ is a smooth map from a neighborhood of
$(0,\x_0)$ in $\R\times(\boW^\pos)^3$ to $\boW$.
Moreover, at $t=0$, we have
$$\wtM_i(0,\x)=2\pi\ii\matrix{a_i+b_ip_i & -\lambda^{-1}(2a_ip_i+b_ip_i^2)\\
\lambda b_i & -a_i-b_ip_i}.$$
\end{proposition}
Proof: we first consider the case where the parameter $\x=(a_i,b_i,p_i)_{1\leq i\leq n}$ is constant with
respect to $\lambda$, so $\x\in\C^{3n}$.
Fix $R>\rho$ and let $\A_R$ be the annulus $\frac{1}{R}<|z|<R$ in $\C$.
For $(\mu,\x)$ in a neighborhood of $(0,\x_0)$ in $\C\times\C^{3n}$ and
$\lambda\in\A_R$, we define
$$\xi_{\mu,\x,\lambda}(z)=\matrix{0&\lambda^{-1}dz\\\mu\,\omega_{\x}(z)&0}$$
where $\omega_{\x}$ is defined as in Equation \eqref{eq-omega} except that
$a_i$, $b_i$, $p_i$ are constant complex numbers.
Let $\Phi_{\mu,\x,\lambda}(z)$ be the solution of $d\Phi_{\mu,\x,\lambda}=\Phi_{\mu,\x,\lambda}\xi_{\mu,\x,\lambda}$ in $\wtOmega$
with initial condition $\Phi_{\mu,\x,\lambda}(\widetilde{0})=I_2$.
Let $N_i(\mu,\x,\lambda)=\boM_{\gamma_i}(\Phi_{\mu,\xi,\lambda})$.
By standard O.D.E. theory, each element $N_{i;k\ell}$ of $N_i$ is a holomorphic
function of $\mu$, $\x$ and $\lambda$.
At $\mu=0$, we have
$$\Phi_{0,\x,\lambda}(z)=\matrix{1&\lambda^{-1}z\\0&1}$$
so $N_i(0,\x,\lambda)=I_2$. Hence
$$\wtN_i(\mu,\x,\lambda):=\frac{\lambda}{\mu}\log N_i(\mu,\x,\lambda)$$
extends holomorphically at $\mu=0$ with value
$\wtN_i(0,\x,\lambda)=\lambda\frac{\partial N_i}{\partial\mu}(0,\x,\lambda)$,
and is holomorphic with respect to $(\mu,\x)$ in a neighborhood of $(0,\x_0)$
and $\lambda\in\A_R$ (as a function of several complex variables).
By Proposition \ref{appendix-prop1} in Appendix \ref{appendix-sectionA}:
$$\frac{\partial N_i}{\partial\mu}(0,\x,\lambda)=
\int_{\gamma_i}\Phi_{0,\x,\lambda}\frac{\partial\xi_{\mu,\x,\lambda}}{\partial\mu}\mid_{\mu=0}
\Phi_{0,\x,\lambda}^{-1}.$$
By the Residue Theorem, we obtain
\begin{eqnarray*}
\wtN_i(0,\x,\lambda)&=&2\pi\ii\,\lambda\Res_{p_i}
\matrix{1 & \lambda^{-1}z\\0 &1}
\matrix{0 & 0 \\ \omega_{\x} & 0}\matrix{1 & -\lambda^{-1}z\\ 0 & 1}\\
&=&2\pi\ii\,\Res_{p_i}\matrix{z & -\lambda^{-1} z^2\\
\lambda & -z}\left(\frac{a_i}{(z-p_i)^2}+\frac{b_i}{z-p_i}\right)\\
&=&2\pi\ii\matrix{a_i+b_ip_i & -\lambda^{-1}(2a_ip_i+b_ip_i^2)\\
\lambda b_i & -a_i-b_ip_i}.
\end{eqnarray*}
In the last equation, we have used the following elementary residue computation:
$$\Res_p\frac{z^k}{(z-p)^2}=kp^{k-1}.$$
For $\x\in(\boW^{\pos})^{3n}$, we have
$$\xi_{t,\x}(z,\lambda)=\xi_{t(\lambda-1)^2,\x(\lambda),\lambda}(z)$$
$$\wtM_i(t,\x)(\lambda)=\wtN_i(t(\lambda-1)^2,\x(\lambda),\lambda).$$
Since the linear map
$(t,\x)\mapsto (t(\lambda-1)^2,\x)$ from $\R\times(\boW^\pos)^{3n}$ to $\boW^{3n+1}$ is bounded,
Proposition \ref{proposition-smooth} in Appendix \ref{appendix-smooth} gives that
for $1\leq k,\ell\leq 2$, the map $(t,\x)\mapsto \wtM_{i;k\ell}(t,\x)$ is smooth from a neighborhood of $(0,\x_0)$ in $\R\times(\boW^{\pos})^{3n}$ to $\boW$.
\cqfd
\begin{remark}
Proposition \ref{proposition-Mi} implies in particular that $M_i(t,\x)(1)=0$ and
$\frac{\partial M_i(t,\x)}{\partial \lambda}(1)=0$ as claimed in Point 2 of Section \ref{nnoids-section1}.
\end{remark}
We define the following functions (the $*$ operator is defined in Section \ref{functional-spaces}):
$$\boF_i(t,\x)=\wtM_{i,11}(t,\x)+\wtM_{i,11}(t,\x)^*$$
$$\boG_i(t,\x)=\lambda\left(\wtM_{i,12}(t,\x)+\wtM_{i,21}(t,\x)^*\right).$$
The Regularity and Rescaled Monodromy Problems are equivalent to the following problem:
\begin{equation}
\label{the-equations}
\left\{\begin{array}{ll}
\boF_i(t,\x)=0 & \mbox{ for $1\leq i\leq n-1$}\\
\boG_i(t,\x)=0 & \mbox{ for $1\leq i\leq n-1$}\\
\boH_i(\x)=0 & \mbox{ for $1\leq i\leq 3$}.
\end{array}\right.
\end{equation}
\section{Solving the equations at $t=0$}
\label{nnoids-section5}
\begin{proposition}
\label{proposition-t0}
When $t=0$, Problem \eqref{the-equations} is equivalent to the following
conditions, for $1\leq i\leq n$:
\begin{enumerate}[(i)]
\item $a_i\in\R$ is constant (with respect to $\lambda$).
\item $p_i$ is constant.
\item $\displaystyle b_i=\frac{-2a_i\overline{p_i}}{1+|p_i|^2}$.
\item $\displaystyle\sum_{i=1}^n a_i\pi^{-1}(p_i)=0$
where $\pi^{-1}$ is the inverse stereographic projection from the south pole
given by Equation \eqref{inverse-stereographic}.
\end{enumerate}
\end{proposition}
Proof:
Using Proposition \ref{proposition-Mi}, we obtain:
\begin{equation}
\boF_i(0,\x)=2\pi\ii\left((a_i+b_i p_i)-(a_i+b_ip_i)^*\right)\label{Fi0}
\end{equation}
\begin{equation}
\boG_i(0,\x)=-2\pi\ii(2a_ip_i+b_ip_i^2+b_i^*)\label{Gi0}
\end{equation}
\begin{equation}
\label{sumFi0}
\sum_{i=1}^n\boF_i(0,\x)=2\pi\ii(\boH_2(\x)-\boH_2(\x)^*)
\end{equation}
\begin{equation}
\label{sumGi0}
\sum_{i=1}^n\boG_i(0,\x)=-2\pi\ii(\boH_3(\x)+\boH_1(\x)^*).
\end{equation}
Assume that $t=0$ and let $\x=(a_i,b_i,p_i)_{1\leq i\leq n}$ be a solution of
Problem \eqref{the-equations}.
From Equations \eqref{sumFi0} and \eqref{sumGi0} we infer that $\boF_n(0,\x)=\boG_n(0,\x)=0$.
Let $i\in[1,n]$.
From Equation \eqref{Gi0}, we see that $b_i\in\boW^{\neg}\cap\boW^{\pos}=\boW^0$
hence $b_i$ is constant.
From Equation \eqref{Fi0}, we obtain by the same argument
\begin{equation}
\label{eqci}
a_i+b_ip_i=c_i
\end{equation}
where $c_i$ is a real constant.
Eliminating $a_i$ from Equations \eqref{Gi0} and \eqref{eqci}, we obtain
$$-b_ip_i^2+2c_i p_i+\overline{b_i}=0.$$
If $b_i\neq 0$ then $p_i$ can take only two values. Being a holomorphic function of $\lambda$, $p_i$ must be constant. By Equation \eqref{eqci}, $a_i$ is constant.
(If $b_i=0$, then $a_i=c_i$ and since $a_i\neq 0$,
Equation \eqref{Gi0} implies that $p_i=0$ is constant.)
Multipliying \eqref{Gi0} by
$\overline{p_i}$ we obtain
\begin{equation}
\label{eqtruc}
(2a_i+b_ip_i)|p_i|^2+\overline{b_ip_i}=0.
\end{equation}
Taking the imaginary part and using $\Im(b_ip_i)=-\Im(a_i)$, we obtain
$$\Im(a_i)(|p_i|^2+1)=0.$$ 
Hence $a_i\in\R$ and so $b_ip_i\in\R$.
Equation \eqref{eqtruc} gives
$$b_i=\frac{-2a_i\overline{p_i}}{1+|p_i|^2}.$$
With this value for $b_i$ we obtain
\begin{equation}
\label{eq-H1}
\boH_1=\sum_{i=1}^n \frac{-2a_i\,\overline{p_i}}{1+|p_i|^2}\end{equation}
\begin{equation}
\label{eq-H2}
\boH_2=\sum_{i=1}^n \frac{a_i(1-|p_i|^2)}{1+|p_i|^2}\end{equation}
\begin{equation}
\label{eq-H3}
\boH_3=\sum_{i=1}^n \frac{2a_i\,p_i}{1+|p_i|^2}.\end{equation}
so the equations $\boH_2=0$ and $\boH_3=0$ give Point (iv).\medskip

Conversely, assume that $a_i$, $b_i$ and $p_i$ satisfy Points (i) to (iv) of Proposition \ref{proposition-t0}.
Proposition \ref{proposition-Mi} gives
$$\wtM_i(0,\x)=2\pi\ii\frac{a_i}{(1+|p_i|^2)}\matrix{1-|p_i|^2&-2\lambda^{-1}p_i\\
-2\lambda\overline{p_i}&|p_i|^2-1}\in\Lambda \su(2).$$
Equations \eqref{eq-H1}, \eqref{eq-H2} and \eqref{eq-H3} imply that $\boH_1=\boH_2=\boH_3=0$.
\cqfd
\section{Solving the equations using the Implicit Function Theorem}
\label{nnoids-section6}
\def\da{{da}}
\def\db{{db}}
\def\dpp{{dp}}
We shall apply the Implicit Function Theorem at the point $(t,\x)=(0,\x_0)$ where
$\x_0$ denotes the following value of the parameters:
\begin{equation}
\label{central-value}
a_i=\tau_i,\quad b_i=\frac{-2a_i\overline{\pi_i}}{1+|\pi_i|^2}
\quad\mbox{ and }\quad p_i=\pi_i
\qquad\mbox{ for $i\in[1,n]$.}
\end{equation}
According to Proposition \ref{proposition-t0} and by the balancing condition of
Theorem \ref{theorem1}, Problem \eqref{the-equations} is solved at $(0,\x_0)$. 
\begin{proposition}
\label{proposition-implicit}
For $t$ in a neighborhood of $0$, there exists a unique smooth map
$t\mapsto \x(t)=(a_{i,t},b_{i,t},p_{i,t})_{1\leq i\leq n}$ with value in $(\boW^{\pos})^{3n}$ 
such that  $\x(0)=\x_0$, Problem \eqref{the-equations} is solved at $(t,\x(t))$
and the following normalisation holds:
\begin{equation}
\label{normalisation}
\forall i\in[1,n-1],\quad
\Re(a_{i,t}^0)=\tau_i\quad\mbox{ and }\quad
p_{i,t}^0=\pi_i.
\end{equation}
\end{proposition}
Proof:
We compute the partial differential of Equations \eqref{Fi0} and \eqref{Gi0}
with respect to $\x$:
$$d_\x\boF_i(0,\x_0)=
2\pi\ii\left(\da_i+p_i\db_i+b_i\dpp_i\right)
-2\pi\ii\left(\da_i+p_i\db_i+b_i\dpp_i\right)^*$$
$$d_\x\boG_i(0,\x_0)=-2\pi\ii\left(
2p_i\da_i+p_i^2\db_i+2(a_i+b_ip_i)\dpp_i+\db_i^*\right).$$
Here, $a_i$, $b_i$, $p_i$
are given by Equation \eqref{central-value} so are constant with respect to $\lambda$.
Projecting on $\boW^{\pospos}$, $\boW^0$ and $\boW^{\negneg}$ we obtain:
$$d_\x\boF_i(0,\x_0)^+=2\pi\ii\left(\da_i^{\pospos}+p_i\db_i^{\pospos}+b_i\dpp_i^{\pospos}\right)$$
$$d_\x\boF_i(0,\x_0)^0=-4\pi\,\Im(\da_i^0+p_i\db_i^0+b_i\dpp_i^0)$$
$$d_\x\boG_i(0,\x_0)^+=-2\pi\ii\left(2p_i\da_i^{\pospos}+p_i^2\db_i^{\pospos}+2(a_i+b_ip_i)\dpp_i^{\pospos}\right)$$
$$d_\x\boG_i(0,\x_0)^0=-2\pi\ii\left(2p_i\da_i^0+p_i^2\db_i^0+2(a_i+b_ip_i)\dpp_i^0+\overline{\db_i^0}\right)$$
$$d_\x\boG_i(0,\x_0)^-=-2\pi\ii(\db_i^{\pospos})^*$$
$$(d_\x\boG_i(0,\x_0)^-)^*=2\pi\ii\db_i^{\pospos}.$$
By definition, we have $\boF_i(t,\x)=\boF_i^*(t,\x)$ for all $(t,\x)$, so $\boF_i^0(t,\x)\in\R$ and $\boF_i^-(t,\x)=0\Leftrightarrow \boF_i^+(t,\x)=0$.
We restrict the parameter $\x$ to the subspace defined by 
Equation \eqref{normalisation}.
\begin{claim}
\label{claim4}
\begin{enumerate}[1.]
\item For $1\leq i\leq n-1$,
the partial differential of $(\boF_i^+,\boG_i^+,(\boG_i^-)^*)$
with respect to $(a_i^+,b_i^+,p_i^+)$ is an automorphism
of $(\boW^\pospos)^3$.
\item For $1\leq i\leq n-1$,
the partial differential of
$(\boF_i^0,\boG_i^0)$ with respect to $(\Im(a_i^0),b_i^0)$ is an automorphism
of $\R\times\C$.
\item The partial differential
of $(\boH_1,\boH_2,\boH_3)$ with respect to $(a_n,b_n,p_n)$ is
an automorphism of $(\boW^{\pos})^3$.
\end{enumerate}
\end{claim}
Proof:
\begin{enumerate}[1.]
\item We can write in matrix form
$$\left(\begin{array}{c}d_{\x}\boF_i^{\pospos}\\ d_{\x}\boG_i^{\pospos}\\ (d_{\x}\boG_i^{\negneg})^*\end{array}\right)
=2\pi\ii
\left(\begin{array}{ccc}1 & p_i & b_i\\ -2p_i&-p_i^2&-2(a_i+b_ip_i)\\
0&1&0\end{array}\right)\left(\begin{array}{c}\da_i^{\pospos}\\ \db_i^{\pospos}\\ \dpp_i^{\pospos}\end{array}\right).$$
This constant matrix has determinant $2a_i$ so is invertible.
(Observe that this operator admits a matrix with respect to
a decomposition of the space as a finite product of Banach spaces. It is clear that
if the matrix is invertible, the operator is an automorphism.)
\item
At $t=0$ and for fixed value of $\Re(a_i)$ and $p_i$,
$(\boF_i(0,\x)^0,\boG_i(0,\x)^0)$ is an affine function of
$(\Im(a_i^0),b_i^0)$. The proof of Proposition \ref{proposition-t0} shows that this function is injective,
so its linear part  is an automorphism of $\R\times\C$.
\item 
The partial differential
of $(\boH_1,\boH_2,\boH_3)$ with respect to $\y=(a_n,b_n,p_n)$ can be written in matrix form as
$$
\left(\begin{array}{c}
d_\y\boH_1\\
d_\y\boH_2\\
d_\y\boH_3\end{array}\right)=
\left(\begin{array}{ccc}
0&1&0\\1&p_n&b_n\\2p_n&p_n^2&2(a_n+b_np_n)\end{array}\right)
\left(\begin{array}{c}\da_n\\ \db_n\\ \dpp_n\end{array}\right).$$
This matrix has determinant $-2a_n\neq 0$.\cqfd
\end{enumerate}
Returning to the proof of Proposition \ref{proposition-implicit}, the partial differential of
$$\left[(\boF_i^+,\boG_i^+,(\boG_i^-)^*,\boF_i^0,\boG_i^0)_{1\leq i\leq n-1},
\boH_1,\boH_2,\boH_3\right]$$
with respect to
$$\left[(a_i^+,b_i^+,p_i^+,\Im(a_i^0),b_i^0)_{1\leq i\leq n-1},a_n,b_n,p_n\right]$$
has lower triangular block form, with automorphisms on the diagonal so is
an automorphism.
(Here I am not talking about a matrix block decomposition but about the decomposition
of an operator with respect to a product of Banach spaces.)
Proposition \ref{proposition-implicit} follows from the Implicit Function Theorem.
\cqfd
%\begin{remark} We have fixed the value of $3n-3$ real parameters, namely $\Re(a_i^0)$
%and $p_i^0$ for $1\leq i\leq n-1$. These parameters
%correspond to deformations of the configuration $(u_i,\tau_i)_{1\leq i\leq n}$
%satisfying the balancing condition. So we have constructed a smooth $3n-3$ dimensional
%family of $n$-noids, or $3n-6$ if we consider surfaces modulo rigid motions.
%This is the dimension of the moduli space of CMC-1 surfaces with $n$ Delaunay ends
%\cite{kusner-mazzeo-pollack}.
%(It looks like we forgot the parameter $t$ in this count, but scaling all weights $\tau_i$ by 
%$k$ gives the same family of surfaces with $t$ replaced by $t/k$.)
%\end{remark}
\section{Delaunay ends}
\label{section-delaunay-ends}
From now on, we assume that the parameter $\x=\x(t)$ is given as a function of $t$ by Proposition \ref{proposition-implicit}.
We denote $\xi_t=\xi_{t,\x(t)}$, $\Phi_t=\Phi_{t,\x(t)}$
and $f_t=\Sym(\Uni(\Phi_t))$ the immersion obtained by the DPW method on $\wtOmega$. Since the Monodromy Problem is solved, $f_t$ descends to a well defined immersion
in $\Omega$, still denoted the same.
In this section we prove:
\begin{proposition}
\label{proposition-delaunay}
\begin{enumerate}[1.]
 For $t\neq 0$ small enough:
\item $f_t$ extends analytically to
\begin{equation}
\label{equation-Sigmat}
\Sigma_t:=\C\cup\{\infty\}\setminus\{p_{1,t}(0),\cdots,p_{n,t}(0)\}.
\end{equation}
\item For $i\in[1,n]$, $a_{i,t}$ is a real constant (with respect to $\lambda$).
\item For $i\in[1,n]$, $f_t$ has a Delaunay end of weight $8\pi ta_{i,t}$ at
$p_{i,t}(0)$.
\end{enumerate}
\end{proposition}
Proof:
\begin{enumerate}[1.]
\item Since the Regularity Problem is solved, the immersion $f_t$ extends analytically at $\infty$.
With the notations of Theorem \ref{theorem-change-coordinate},
we consider the change of variable
$$z=\psi_{i,t,\lambda}(w)=p_{i,t}(\lambda)+w$$
and the following domains:
$$U=\{w\in\C:2\varepsilon<|w|<4\varepsilon\}$$
$$V_i=\{z\in\C:\varepsilon<|z-\pi_i|<8\varepsilon\}\subset\Omega$$
For $t$ small enough, we have
$\psi_{i,t,\lambda}(U)\subset V_i$ for all $\lambda\in\D_{\rho}$.
Let $\wtU$ and $\wtV_i$ be the universal covers of $U$ and $V_i$.
(For $\wtV_i$, we may take an arbitrary component of $\wtOmega\cap p^{-1}(V_i)$
where $p:\wtOmega\to\Omega$ is the universal cover.)
Lift $\psi_{i,t,\lambda}$ to $\wtpsi_{i,t,\lambda}:\wtU\to\wtV_i$ and define
for $w\in\wtU$:
$$\whPhi_{i,t}(w,\lambda)=\Phi_t(\wtpsi_{i,t,\lambda}(w),\lambda).$$
By Corollary \ref{corollary-change-coordinate}, 
$\whf_{i,t}=\Sym(\Uni(\whPhi_{i,t}))$ descends to a well defined immersion on $U$
and
\begin{equation}
\label{eq-wtf}
\forall w\in U,\quad
\whf_{i,t}(w)=f_t(\psi_{i,t,0}(w))=f_t(p_{i,t}(0)+w).\end{equation}
Now $\whPhi_{i,t}$ solves $d\whPhi_{i,t}=\whPhi_{i,t}\whxi_{i,t}$ with
$$\whxi_{i,t}=\psi_{i,t,\lambda}^*\xi_t=\matrix{0& \lambda^{-1}dw\\t(\lambda-1)^2\whomega_{i,t}&0}
\quad\mbox{ and }\quad \whomega_{i,t}=\psi_{i,t,\lambda}^*\omega_{\x(t)}.$$
Since the only pole of $\whomega_{i,t}$ in $D(0,4\varepsilon)$ is at $w=0$, the DPW method shows that $\whf_{i,t}$ extends analytically to $D^*(0,4\varepsilon)$.
We may extend $f_t$ analytically by requesting that Equation \eqref{eq-wtf} holds true in $D^*(0,4\varepsilon)$.
Doing this for $i\in[1,n]$, we have extended $f_t$ analytically to $\Sigma_t$.
\item Next we gauge $\whxi_{i,t}$ so that it has a simple pole at $0$, yielding a fuchsian system
of differential equations. Consider the gauge
$$G(w,\lambda)=\matrix{\frac{\sqrt{w}}{k}&0\\-\frac{\lambda}{2k\sqrt{w}}&\frac{k}{\sqrt{w}}}.$$
Here we can take $k=1$, but in the next point we will take another value of $k$ so
we do the computation for general values of $k\neq 0$.
The gauged potential is
$$\wcxi_{i,t}:=\whxi_{i,t}\cdot G=\matrix{0&\lambda^{-1}k^2\frac{dw}{w}\\
\frac{t(\lambda-1)^2}{k^2}w\,\whomega_{i,t}+\frac{\lambda}{4k^2}\frac{dw}{w}&0}.$$
It has a simple pole at $0$ with residue
\begin{equation}
\label{eq-Ait}
A_{i,t}(\lambda):=\matrix{0&\lambda^{-1}k^2\\\frac{t(\lambda-1)^2}{k^2}a_{i,t}(\lambda)
+\frac{\lambda}{4k^2}&0}.
\end{equation}
The eigenvalues of $A_{i,t}(\lambda)$ are $\pm\Lambda_{i,t}(\lambda)$ with
$$\Lambda_{i,t}(\lambda)^2=t\lambda^{-1}(\lambda-1)^2a_{i,t}(\lambda)+\frac{1}{4}.$$
Fix $\lambda\in\S^1\setminus\{1\}$. If $t\neq 0$ is small enough, then $\Lambda_{i,t}(\lambda)\not\in
\Z/2$ so the corresponding fuchsian system is non-resonant.
Hence $\wcPhi_{i,t}:=\whPhi_{i,t} G$
has the following standard $z^AP$ form in the universal cover of $D(0,\varepsilon)^*$ (see Proposition 11.2 in \cite{taylor})
$$\wcPhi_{i,t}(w,\lambda)=V(\lambda)\exp(A_{i,t}(\lambda)\log w)P(w,\lambda)$$
where $V\in\Lambda SL(2,\C)$ and $P(w,\lambda)$ descends to a well defined holomorphic function of $w$ in
$D(0,\varepsilon)$ with $P(0,\lambda)=I_2$.
Hence
$$\boM_{C(0,\varepsilon)}(\wcPhi_{i,t})(\lambda)=V(\lambda)\exp(2\pi\ii A_{i,t}(\lambda))V(\lambda)^{-1}$$
so the eigenvalues of $\boM_{\gamma_i}(\Phi_t)$ are
$\exp(\pm 2\pi\ii\Lambda_{i,t}(\lambda))$.
Since the Monodromy Problem is solved, the eigenvalues are unitary complex numbers,
so $\Lambda_{i,t}(\lambda)\in\R$ which implies $a_{i,t}(\lambda)\in\R$.
Since $a_{i,t}$ is holomorphic in $\D_{\rho}$ and $a_{i,t}$ is real on $\S^1\setminus\{1\}$, it is constant.
\item It remains to choose $k$ in the definition of $G$ such that $A_{i,t}(\lambda)$ is a standard Delaunay
residue. Fix a small $t>0$ and
let $(r,s)\in\R^2$ be the solution of the system
\begin{equation}
\label{equation-rs}
\left\{\begin{array}{l}
rs=ta_{i,t}\\
r+s=\frac{1}{2}\\
r>s
\end{array}\right.\end{equation}
For small $t$, $(r,s)$ is close to $(\frac{1}{2},0)$ so
$\sqrt{r+s\lambda}$ is well defined and does not vanish in $\D_{\rho}$.
We take
$$k=\sqrt{r+s\lambda}$$ in the definition of $G$.
Using $r+s=\frac{1}{2}$, we have
\begin{equation}
\label{detArs}
(r+s\lambda)(r\lambda+s)=rs(\lambda-1)^2+\frac{\lambda}{4}=ta_{i,t}(\lambda-1)^2+\frac{\lambda}{4}.\end{equation}
Hence by Equation \eqref{eq-Ait},
\begin{equation}
\label{equation-Ars}
A_{i,t}(\lambda)=
\matrix{0 & r\lambda^{-1}+s\\r\lambda+s & 0}
\end{equation}
which is the residue of the standard Delaunay potential (see Section \ref{section-basic-examples}).
Since the Monodromy Problem is solved, the immersion $\whf_{i,t}$ has
a Delaunay end at $w=0$ of weight $8\pi rs=8\pi ta_{i,t}$ by
Theorem 3.5 in \cite{kilian-rossman-schmitt}.\cqfd
\end{enumerate}
\begin{remark}
We have constructed a smooth family of immersions $f_t$ defined for $t\neq 0$
in a neighborhood of $0$. As $t$ changes sign, the behavior of each end switches
from unduloid to nodoid type. As $t\to 0$, the immersion $f_t$ degenerates into a sphere.
\end{remark}
\section{Geometry of the immersion}
\label{nnoids-section8}
Theorem 3.5 of \cite{kilian-rossman-schmitt} tells us that
for each $t>0$, there exists a Delaunay immersion $f_{i,t}^{\boD}:\C^*\to\R^3$ such that
$$\lim_{w\to 0}||\whf_{i,t}(w)-f_{i,t}^{\boD}(w)||=0.$$
In that paper, $t$ is fixed. The problem is that the limit $w\to 0$ is not uniform with
respect to $t$. Indeed, the fuchsian system is resonant at $t=0$ so the constants in
their estimates are uncontroled as $t\to 0$.
Thomas Raujouan has improved this result in \cite{raujouan}
and was able to obtain a uniform limit under additional assumptions:
\begin{theorem}\cite{raujouan}
\label{theorem-thomas}
Let $\tau$ be a non-zero real number.
Let $\xi_t(z,\lambda)$ be a family of DPW potentials depending on the parameter
$t\geq 0$ and defined for $z$ in a punctured neighborhood of $0$.
Let $\Phi_t(z,\lambda)$ be a solution of $d_z\Phi_t=\Phi_t\xi_t$ and
$f_t(z)=\Sym(\Uni(\Phi_t))$.
Assume that:
\begin{enumerate}[1.]
\item $\xi_t=A_t\frac{dz}{z}+O(t,z^0)$ where $A_t$ is the standard Delaunay residue
given by Equation \eqref{equation-Ars} and $(r,s)\in\R^2$ is the solution of 
$$\left\{\begin{array}{l}
rs=t\tau\\
r+s=\frac{1}{2}\\
r>s\end{array}\right.$$
\item The monodromy of $\Phi_t$ around the origin is in $\Lambda SU(2)$.
\item $\Phi_0(1,\lambda)=\minimatrix{a&\lambda^{-1}b\\\lambda c&d}$
where $a,b,c,d$ are (constant) complex numbers.
\end{enumerate}
Then:
\begin{enumerate}[1.]
\item For all $0<\alpha<1$, there exists uniform positive numbers $\epsilon$, $c$, $T$ and a family of Delaunay immersions
$f_t^{\boD}:\C^*\to\R^3$ with weight $8\pi\tau t$ such that for all $0<t<T$ and
$0<|z|<\epsilon$,
$$||f_t(z)-f_t^{\boD}(z)||\leq c\,t|z|^\alpha.$$
\item Let
$$H(\lambda)=\frac{1}{\sqrt{2}}\matrix{1&-\lambda^{-1}\\\lambda&1}$$
$$Q=\Uni(\Phi_0(1,\cdot)H).$$
The axis of $f_t^{\boD}$ (oriented from the end at $\infty$ to the end at $0$) converges as $t\to 0$ to the line through the point $\Sym(Q)$ spanned by the vector $-\Nor(Q)$
(notations as in Section \ref{background-section3}).
\item If $\tau>0$, there exists $T'<T$ such that the restriction of $f_t$ to $D^*(0,\epsilon)$
is an embedding for all $0<t<T'$. More precisely, $f_t(D^*(0,\epsilon))$
is included in an embedded tubular neighborhood of the Delaunay surface
$f_t^{\boD}(\C^*)$ and the projection on the Delaunay surface is a diffeomorphism
from $f_t(D^*(0,\epsilon))$ onto its image.
\end{enumerate}
\end{theorem}
Observe that Hypothesis 1 implies that $\xi_0=A_0\frac{dz}{z}$ so $\Phi_0$ is actually
defined for $z\in\C^*$, which is why Hypothesis 3 makes sense.
\subsection{Axes of the ends}
\label{section-axes}
We first use Theorem \ref{theorem-thomas} to compute the limit axes of the ends.
\begin{proposition}
\label{proposition-axe}
The axis of the Delaunay end of $f_t$ at $p_{i,t}(0)$ converges as $t\to 0$ to
the half-line through $(0,0,-1)$ spanned by the vector $u_i$.
\end{proposition}
Proof: we continue with the notations of the proof of Proposition \ref{proposition-delaunay}.
We want to apply Theorem \ref{theorem-thomas} to $\wcPhi_{i,t}=\whPhi_{i,t}G$.
First of all, $t\mapsto ta_{i,t}$ is a smooth diffeomorphism in a neighborhood
of $0$ so we may use $|ta_{i,t}|$ as the time parameter in our application
of Theorem \ref{theorem-thomas}, with $\tau=1$ if $ta_{i,t}>0$ and $\tau=-1$ if
$ta_{i,t}<0$.
We have
$$\whPhi_{i,0}(w,\lambda)=\Phi_0(\pi_i+w,\lambda)=\matrix{1&\lambda^{-1}(\pi_i+w)\\0&1}.$$
At $t=0$, the solution of \eqref{equation-rs} is $(r,s)=(\frac{1}{2},0)$ so $k=\frac{1}{\sqrt{2}}$.
This gives
\begin{eqnarray}
\wcPhi_{i,0}(1,\lambda)&=&\whPhi_{i,0}(1,\lambda)G(1,\lambda)\nonumber\\
&=&\matrix{1&\lambda^{-1}(\pi_i+1)\\0&1}\matrix{\sqrt{2}&0\\-\frac{\lambda}{\sqrt{2}}&\frac{1}{\sqrt{2}}}\nonumber\\
&=&\frac{1}{\sqrt{2}}\matrix{1-\pi_i&\lambda^{-1}(\pi_i+1)\\-\lambda&1}.
\label{eq-M}
\end{eqnarray}
Fix some $\alpha\in(0,1)$.
By Theorem \ref{theorem-thomas}, there exists uniform positive numbers $\epsilon$, $c$, $T$ 
and a family of Delaunay immersions
$f_{i,t}^{\boD}$ such that for $0<t<T$ and
$0<|w|<\epsilon$
$$||\whf_{i,t}(w)-f_{i,t}^{\boD}(w)||\leq c\,t|w|^{\alpha}.$$
Using Equation \eqref{eq-M}, the matrix $Q$ in Point 2 of Theorem \ref{theorem-thomas} is given by
$$Q(\lambda)=\Uni\matrix{1&\lambda^{-1}\pi_i\\0&1}.$$
By the results in Section \ref{section-basic-examples} for the standard sphere, the limit axis of $f_{i,t}^{\boD}$ is the line through the point $(0,0,-1)+\pi^{-1}(\pi_i)$
directed by the vector $\pi^{-1}(\pi_i)$. Since $\pi_i=\pi(u_i)$, this proves
Proposition \ref{proposition-axe}.\cqfd
\medskip

\subsection{Embeddedness}
\label{section-embedded}
For ease of notation, we forget from now on the argument $\lambda=0$ so we write
$p_{i,t}=p_{i,t}(0)$.
It will be convenient to translate all immersions by the vector $(0,0,1)$, so
we replace $f_t$ by $f_t+(0,0,1)$ and $f_{i,t}^{\boD}$ by $f_{i,t}^{\boD}+(0,0,1)$.
Thanks to this translation, the standard spherical data (see Section \ref{section-basic-examples})
yields an immersion $f_0:\C\cup\{\infty\}\to\S^2$ equal to the inverse stereographic
projection.
Then $f_t$ converges smoothly on compact subsets of $\C\cup\{\infty\}\setminus
\{\pi_1,\cdots,\pi_n\}$ to $f_0$.
Indeed, on compact subsets of $\C\setminus\{\pi_1,\cdots,\pi_n\}$,
$\xi_t$ converges uniformly to the standard spherical potential
$\xi_0=\minimatrix{0&\lambda^{-1}\\0&0}dz$, with the same initial data;
and in a neighborhood of $\infty$, $\xi_t\cdot G_{\infty}$ converges uniformly
to $\xi_0\cdot G_{\infty}$.
Let $M_t=f_t(\Sigma_t)$, where $\Sigma_t$ is defined by Equation \eqref{equation-Sigmat}.
\begin{proposition}
\label{proposition-embedded}
If all weights $\tau_i$ are positive and the angle between $u_i$ and $u_j$ is greater than $\frac{\pi}{3}$ for all $j\neq i$, then for $t>0$ small enough, $M_t$ is embedded.
\end{proposition}
Proof: we continue with the notations of Section \ref{section-axes}. In particular, $\epsilon$
is the number given by our application of Theorem \ref{theorem-thomas} and is fixed.
\begin{figure}
\begin{center}
\includegraphics[height=5cm]{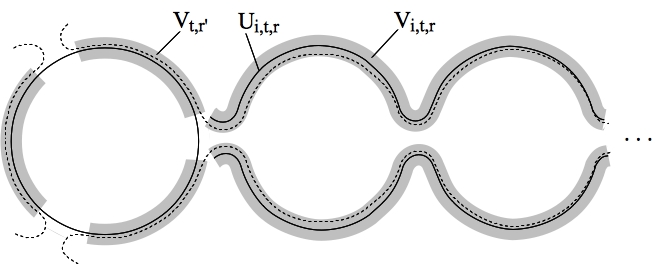}
\end{center}
\caption{The domains $V_{i,t,r}$ and $V_{t,r'}$. On this figure, $n=3$ and $r'>r$ so the domains $V_{i,t,r}$ and $V_{t,r'}$ do not overlap (for clarity). The surface $M_t$ is represented with dots. The Delaunay surface $\boD_{i,t}$ and the sphere $\S^2$
are represented with plain lines.}
\label{figure2}
\end{figure}
\begin{itemize}
\item For $i\in[1,n]$, we denote $\boD_{i,t}$ the Delaunay surface $f_{i,t}^{\boD}(\C^*)$.
By Point 3 of Theorem \ref{theorem-thomas}, for $t$ small enough and for all $r\leq\epsilon$, the projection from $f_t(D^*(p_{i,t},r))$ to $\boD_{i,t}$ is a diffeomorphism onto its image
which we denote $U_{i,t,r}$. In other words, $f_t(D^*(p_{i,t},r))$ is the normal graph
on $U_{i,t,r}\subset \boD_{i,t}$ of a function which, by Point 1 of Theorem \ref{theorem-thomas},
is bounded by $c_1t$ where
$c_1=c\epsilon^{\alpha}$. We denote $V_{i,t,r}=\mbox{Tub}_{c_1t}U_{i,t,r}$ the
tubular neighborhood of radius $c_1t$ of $U_{i,t,r}$. We have
$f_t(D^*(p_{i,t},r))\subset V_{i,t,r}$ (see Figure \ref{figure2}).
\item For $r>0$, let
$$\Omega_{t,r}=\{z\in\C:\forall i\in[1,n],\;|z-p_{i,t}|>r\}\cup\{\infty\}.$$
For $t>0$ small enough, $f_t(\Omega_{t,\frac{\epsilon}{8}})$ is the normal graph,
on a domain of the sphere $\S^2$, of a function which is bounded by $c_2 t$ for some
constant $c_2$. For $r\in[\frac{\epsilon}{8},\epsilon]$, we denote $U_{t,r}$ the
projection of $f_t(\Omega_{t,r})$ on the sphere and $V_{t,r}=\mbox{Tub}_{c_2t}U_{t,r}$
the tubular neighborhood of radius $c_2t$ of $U_{t,r}$. We have
$f_t(\Omega_{t,r})\subset V_{t,r}$.
\end{itemize}
For $p\in\R^3$, we denote $\S^2(p)$ the unit sphere centered at $p$.
\begin{claim}
For $r\in[\frac{\epsilon}{8},\epsilon]$:
\begin{equation}
\label{eq-Lr}
\lim_{t\to 0}\overline{V}_{t,r}=\S^2(0)\setminus\bigcup_{i=1}^n f_0(D(\pi_i,r))=:L_r
\end{equation}
\begin{equation}
\label{eq-Lir}\forall i\in[1,n],\quad \lim_{t\to 0}\overline{V}_{i,t,r}=f_0(\overline{D}(\pi_i,r))\cup
\bigcup_{k\geq 1}\S^2(2ku_i)=:L_{i,r}.
\end{equation}
In Equation \eqref{eq-Lr}, the limit is for the Hausdorff distance $d_{\boH}$,
and in Equation \eqref{eq-Lir} for the Hausdorff distance on compact subsets of
$\R^3$: $\lim_{t\to 0} A_t=L$ means that for all compact $K\subset\R^3$,
$\lim_{t\to 0} d_{\boH}( A_t\cap K,L\cap K)=0$.
\end{claim}
Proof: 
\begin{enumerate}[1.]
\item 
The Hausdorff distance between $\overline{V}_{t,r}$ and
$f_t(\overline{\Omega}_{t,r})$ is at most $2c_2 t$, so
$$\lim_{t\to 0}\overline{V}_{t,r}=\lim_{t\to 0} f_t(\overline{\Omega}_{t,r})
=f_0(\C\cup\{\infty\}\setminus\bigcup_{i=1}^n D(\pi_i,r))=L_r.$$
\item It is known that as $t\to 0$, the family of Delaunay surfaces $\boD_{i,t}$ converges, up to rigid motions, to an infinite chain of unit-spheres. Since the axis of $\boD_{i,t}$ converges as $t\to 0$ to the line spanned by $u_i$, the family $\boD_{i,t}$ converges,
up to translation along this line, to a chain of unit-spheres with centers on this line.
We have
$$\lim_{t\to 0}f_{i,t}^{\boD}( \overline{D}(0,r)\setminus D(0,\smallfrac{r}{2}))
=\lim_{t\to 0}f_t( \overline{D}(p_{i,t},r)\setminus D(p_{i,t},\smallfrac{r}{2}))
=f_0(\overline{D}(\pi_i,r)\setminus D(\pi_i,\smallfrac{r}{2}))
\subset \S^2(0).$$
Hence $\boD_{i,t}$ converges to the union of the spheres $\S^2(2k u_i)$ for $k\in\Z$
and
$$\lim_{t\to 0}\overline{V}_{i,t,r}
=\lim_{t\to 0}\overline{f_t( D^*(p_{i,t},r))}
=\lim_{t\to 0}\overline{f_{i,t}^{\boD}(D^*(0,r))}=f_0(\overline{D}(\pi_i,r))\cup
\bigcup_{k\geq 1}\S^2(2ku_i)=L_{i,r}.$$
\cqfd
\end{enumerate}
\begin{claim}
There exists a positive $T'<T$ such that for $0<t<T'$ and $r\in[\frac{\epsilon}{4},\epsilon]$:
\begin{equation}
\label{equation-ViVj}
\forall j\neq i,\quad \overline{V}_{i,t,r}\cap \overline{V}_{j,t,r}=\emptyset
\end{equation}
\begin{equation}
\label{equation-ViV}\forall i,\quad \overline{V}_{i,t,\frac{r}{2}}\cap \overline{V}_{t,r}=\emptyset
\end{equation}
\end{claim}
Proof: The angle hypothesis of Proposition \ref{proposition-embedded} ensures that
the spheres $\S^2(2ku_i)$ and $\S^2(2\ell u_j)$ for
$j\neq i$ and $k,\ell\geq 1$ do not intersect.
Hence the limit sets $L_{i,r}$
for $i\in[1,n]$ are disjoint. Let $\delta>0$ be the smallest distance between $L_{i,r}$
and $L_{j,r}$ for $j\neq i$. Let $K=\overline{B}(0,4)$.
There exists a positive $T'\leq T$ such that for $t<T'$ and $i\in[1,n]$,
$$d_{\boH}(\overline{V}_{i,t,r}\cap K, L_{i,r}\cap K)<\frac{\delta}{2}$$
Hence for $t<T'$ and $j\neq i$, $\overline{V}_{i,t,r}$ and $\overline{V}_{j,t,r}$ do not intersect
inside $K$. Since it is clear that they don't intersect outside $K$,
Equation \eqref{equation-ViVj} follows.
Equation \eqref{equation-ViV} is proved in the same way, observing that
the limit sets $L_{i,\frac{r}{2}}$ and $L_{r}$ are disjoint. 
\cqfd
\begin{claim}
\label{claim-embedded}
For $0<t<T'$ and $i\in[1,n]$:
$$
f_t^{-1}(V_{i,t,\frac{\epsilon}{2}})=D^*(p_{i,t},\smallfrac{\epsilon}{2})
$$
$$
f_t^{-1}(V_{t,\frac{\epsilon}{4}})=\Omega_{t,\frac{\epsilon}{4}}
$$
\end{claim}
Proof:
\begin{enumerate}[1.]
\item Let $z\in f_t^{-1}(V_{i,t,\frac{\epsilon}{2}})$. By
Equation \eqref{equation-ViV} with $r=\epsilon$, $f_t(z)\not\in \overline{V}_{t,\epsilon}$
hence $z\not\in\overline{\Omega}_{t,\epsilon}$.
So there exists $j\in[1,n]$ such that $z\in D^*(p_{j,t},\epsilon)$.
Equation \eqref{equation-ViVj} with $r=\epsilon$ yields $j=i$.
Since the projection from $f_{i,t}(D^*(p_{i,t},\epsilon))$ to $U_{i,t,\epsilon}$ is a diffeomorphism, $z\in D^*(p_{i,t},\frac{\epsilon}{2})$.
\item
Let $z\in f_t^{-1}(V_{t,\frac{\epsilon}{4}})$.
Equation \eqref{equation-ViV} with $r=\frac{\epsilon}{4}$ implies that for all $i\in[1,n]$,
$f_t(z)\not\in \overline{V}_{i,t,\frac{\epsilon}{8}}$ so
$z\not\in\overline{D}(p_{i,t},\frac{\epsilon}{8})$. Hence $z\in\Omega_{t,\frac{\epsilon}{8}}$.
Since the projection from $f_t(\Omega_{t,\frac{\epsilon}{8}})$ to $U_{t,\frac{\epsilon}{8}}$ is a diffeomorphism, 
$z\in\Omega_{t,\frac{\epsilon}{4}}$.\cqfd
\end{enumerate}
Claim \ref{claim-embedded} implies that $M_t$ is embedded. Indeed,
the open sets $V_{i,t,\frac{\epsilon}{2}}$ for $i\in[1,n]$ and $V_{t,\frac{\epsilon}{4}}$ cover 
$M_t$ and the intersection of $M_t$ with each of these sets is a submanifold,
so $M_t$ is a submanifold of $\R^3$.
\cqfd

\subsection{Alexandrov-embeddedness}
We recall from \cite{karsten-kusner-sullivan} the definition of Alexandrov-embeddedness in the non-compact case:
\begin{definition} A CMC surface $M$ of finite topology is Alexandrov-embedded if M is properly immersed, if each end of $M$ is embedded, and if there exists a compact 3-manifold $W$ with
boundary $\partial W = \Sigma$ and a proper immersion $F : W \setminus\{q_1,\cdots, q_n\} \to\R^3$ whose boundary restriction $f : \Sigma\setminus\{q_1,\cdots,q_n\}\to\R^3$ parametrizes $M$. Moreover, we require that the mean-curvature
normal of $M$ points into $W$.
\end{definition}
\begin{proposition}
If all weights $\tau_i$ are positive then for $t>0$ small enough, $M_t$ is Alexandrov-embedded.
\end{proposition}
Proof: the idea is to construct an abstract flat $3$-manifold $N_t$ in which the $n$ half
Delaunay surfaces do not intersect, so $M_t$ lifts to an embedded surface in $N_t$.
\medskip

We continue with the notations of Section \ref{section-embedded}.
Let $\boC_{i,t}$ be the solid cylinder in $\R^3$ bounded by the Delaunay surface $\boD_{i,t}$. For positive $\delta$, let
$$Q_{i,t,\delta}=\{x\in\R^3:d(x,\boC_{i,t})<c_1 t\mbox{ and } \langle x,u_i\rangle>1-\delta\}.$$
Then
$$\lim_{t\to 0} \overline{Q}_{i,t,\delta}=\{x\in \overline{B}(0,1): \langle x,u_i\rangle\geq 1-\delta\}\cup
\bigcup_{k\geq 1}\overline{B}(2ku_i,1)=: L_{i,\delta}.$$
We have
$$\lim_{\delta\to 0}L_{i,\delta}\cap\overline{B}(0,1+\delta)=\{u_i\}$$
so we may choose $\delta>0$ small enough such that the sets
$L_{i,\delta}\cap \overline{B}(0,1+\delta)$ for $i\in[1,n]$ are
disjoint. Then for $t$ small enough and $j\neq i$,
$$Q_{j,t,\delta}\cap Q_{i,t,\delta}\cap B(0,1+\delta)=\emptyset.$$
Let $N_t$ be the flat $3$-manifold obtained as the disjoint union of $B(0,1+\delta)$
and $Q_{i,t,\delta}$ for $i\in[1,n]$, gluing $B(0,1+\delta)$ with $Q_{i,t,\delta}$
(with the identity map) where they intersect in $\R^3$.
Let $\psi_t:N_t\to\R^3$ be the canonical projection. Then $\psi_t$ is an immersion and
$N_t$ is called an immersed domain.
The point here is that the domains $Q_{i,t,\delta}$ for $i\in[1,n]$ are disjoint in $N_t$,
even if they may intersect in $\R^3$.
\medskip

We have for $t$ small enough
$$f_t(\Omega_{t,\frac{\epsilon}{2}})\subset B(0,1+\delta)$$
and taking $\epsilon$ smaller if necessary,
$$f_t(D^*(p_{i,t},\epsilon))\subset Q_{i,t,\delta}.$$
Hence we may lift the immersion $f_t$ to $\wtf_t:\Sigma_t\to N_t$ so that
$\psi_t\circ \wtf_t=f_t$. Let $\wtM_t=\wtf_t(\Sigma_t)$.
Since the domains $Q_{i,t,\delta}$ for $i\in[1,n]$ are disjoint in $N_t$, the proof of
Proposition \ref{proposition-embedded} gives that $\wtM_t$ is embedded in $N_t$.
Let $W_t$ be the domain in $N_t$ bounded by $\wtM_t$.
Since the Delaunay surfaces are Alexandrov-embedded
(in the positive weight case), we may compactify $W_t$ by adding one point at
infinity per Delaunay end of $\wtM_t$. We take $F_t$ to be the restriction of
$\psi_t$ to $W_t$. This proves that $M_t$ is Alexandrov embedded.
\cqfd
\subsection{Hopf differential and umbilics}
\label{umbilics}
Umbilics are points where the two principal curvatures are equal.
On a CMC-1 surface, they are the zeros of the Hopf quadratic differential.
Writing the potential $\xi$ as in Equation \eqref{equation-xi},
the Hopf differential is given by
$$Q(z)=-2\beta(z,0)\gamma(z,0).$$
In our case, the Hopf differential is equal to $-2t\omega_t(z,0)\,dz$.
Since $\omega_t$ has a double zero at $\infty$, the Hopf differential is holomorphic
at $\infty$ and has $n$ double poles at the punctures.
Being a quadratic
differential on the Riemann sphere, it has $2n-4$ zeros, so there are $2n-4$
umbilics (counting multiplicity).
The umbilics converge as $t\to 0$ to the zeros of $\omega_0\,dz$, which, by 
Equation \eqref{central-value}, is given by
$$\omega_0\,dz=\sum_{i=1}^n\left(\frac{\tau_i}{(z-\pi_i)^2}
-\frac{2\tau_i\overline{\pi_i}}{(1+|\pi_i|^2)(z-\pi_i)}\right)dz^2.$$
So the limit position of the umbilics can in principle be computed by solving polynomial
equations.
\appendix
\section{derivative of the monodromy}
\label{appendix-sectionA}
The following proposition is adapted from Proposition 9 in \cite{traizet}.
\begin{proposition}
\label{appendix-prop1}
 Let $\xi_t$ be a $C^1$ family of matrix-valued 1-forms on a Riemann surface $\Sigma$.
Let $\wtSigma$ be the universal cover of $\Sigma$.
Fix a point $z_0$ in $\Sigma$ and let $\wtz_0$ be a lift of $z_0$ to $\wtSigma$.
Let $\Phi_t$ be a family of solutions of $d\Phi_t=\Phi_t\xi_t$ on $\wtSigma$,
such that $\Phi_t(\wtz_0)$ does not depend on $t$. 
Let $\gamma\in\pi_1(\Sigma,z_0)$ and let $M(t)$
be the monodromy of $\Phi_t$ with respect to $\gamma$.
Let $\wtgamma$ be the lift of $\gamma$ to $\wtSigma$ such that
$\wtgamma(0)=\wtz_0$.
Then for all $t$,
$$M'(t)=\int_{\wtgamma}\Phi_{t}\frac{\partial\xi_t}{\partial t}\Phi_{t}^{-1}
\times M(t).$$
\end{proposition}
Proof: since $\Phi_t(\wtz_0)$ is constant and $\xi_t$ depends $C^1$ on $t$, $\Phi_t$ depends $C^1$ on $t$.
Let $\Psi_t=\frac{\partial\Phi_t}{\partial t}$.
By differentiation of the Cauchy Problem satisfied by $\Phi_t$ with respect to $t$, 
we obtain that $\Psi_t$ satisfies the following Cauchy Problem on $\wtSigma$:
$$\left\{\begin{array}{l}
d\Psi_t=\Psi_t\xi_{t}+\Phi_{t}\frac{\partial\xi_t}{\partial t}\\
\Psi_t(\wtz_0)=0.\end{array}\right.$$
Following the method of variation of constants , the function $U_t=\Psi_t\Phi_t^{-1}$
satisfies:
$$\left\{\begin{array}{l}
dU_t=\Phi_t\frac{\partial\xi_t}{\partial t}\Phi_t^{-1}\\
U_t(\wtz_0)=0.\end{array}\right.$$
Hence (writing $\wtgamma(1)$ for the endpoint of $\wtgamma$)
$$U_t(\wtgamma(1))=\int_{\wtgamma}\Phi_t\frac{\partial\xi_t}{\partial t}\Phi_t^{-1}.$$
We have by definition
$$M(t)=\Phi_t(\wtgamma(1))\Phi_t(\wtz_0)^{-1}.$$
Hence since $\Phi_t(\wtz_0)$ is constant:
\begin{eqnarray*}
M'(t)&=&\Psi_t(\wtgamma(1))\Phi_t(\wtz_0)^{-1}\\
&=&U_t(\wtgamma(1))\Phi_t(\wtgamma(1))\Phi_t(\wtz_0)^{-1}\\
&=&\int_{\gamma}
\Phi_{t}\frac{\partial\xi_t}{\partial t}\Phi_t^{-1} M(t).
\end{eqnarray*}
\cqfd
\section{Smoothness of maps between Banach spaces}
\label{appendix-smooth}
The following proposition is useful to prove that the maps considered in
this paper are smooth maps between Banach spaces.
The Banach algebra $\boW$ is defined in Section \ref{functional-spaces}.
For $R>1$, we denote $\A_R$ the annulus $\frac{1}{R}<|\lambda|<R$ in $\C$.
For ${\bf a}=(a_1,\cdots,a_n)\in\C^n$ and ${\bf r}=(r_1,\cdots,r_n)\in(0,\infty)^n$,
we denote $D({\bf a},{\bf r})$ the polydisk $\prod_{i=1}^n D(a_i,r_i)$ in $\C^n$.
\begin{proposition}
\label{proposition-smooth}
Let $R>\rho$ and $f:\A_R\times D({\bf a},{\bf r})\to\C$ be a holomorphic function of 
$(n+1)$ variables $(\lambda,z_1,\cdots,z_n)$.
Let
$$B({\bf a},{\bf r})=\{(u_1,\cdots,u_n)\in \boW^n:\forall i\in[1,n], ||u_i-a_i||<r_i\}$$
where we identify $a_i$ with a constant function in $\boW$.
Define for $(u_1,\cdots,u_n)\in B({\bf a},{\bf r})$:
$$F(u_1,\cdots,u_n)(\lambda)=f(\lambda,u_1(\lambda),\cdots,u_n(\lambda)).$$
Then
$F:B({\bf a},{\bf r})\subset\boW^n\to\boW$ is of class $C^{\infty}$.
\end{proposition}
Proof: we expand $f$ in Laurent series with respect to $\lambda$
and power series with respect to $z_1,\cdots,z_n$:
$$f(\lambda,z_1,\cdots,z_n)=\sum_{k\in\Z}\sum_{i_1,\cdots,i_n}
c_{k\,i_1\cdots i_n}\lambda^k(z_1-a_1)^{i_1}\cdots (z_n-a_n)^{i_n}.$$
For any ${\bf r}'<{\bf r}$ (in the sense $r'_i<r_i$ for all $i$), the series $f(\rho^{\pm 1},a_1+r'_1,\cdots,a_n+r'_n)$ converges absolutely so
\begin{equation}
\label{convergence-norme}
\sum_{k\in\Z}\sum_{i_1,\cdots,i_n}
|c_{k\,i_1\cdots i_n}|\rho^{|k|}(r'_1)^{i_1}\cdots (r'_n)^{i_n}<\infty.
\end{equation}
Let $v\in\boW$ be the function defined by $v(\lambda)=\lambda$. Then formally:
\begin{equation}
\label{formal-series}F(u_1,\cdots,u_n)=\sum_{k\in\Z}\sum_{i_1,\cdots,i_n}
c_{k\,i_1\cdots i_n}v^k(u_1-a_1)^{i_1}\cdots (u_n-a_n)^{i_n}.
\end{equation}
Since $||v^k||=\rho^{|k|}$ and $\boW$ is a Banach algebra, for $(u_1,\cdots,u_n)\in B({\bf a},{\bf r}')$, we have
$$||c_{k\,i_1\cdots i_n}v^k(u_1-a_1)^{i_1}\cdots (u_n-a_n)^{i_n}||\leq |c_{k\,i_1\cdots i_n}|\rho^{|k|}(r'_1)^{i_1}\cdots(r'_n)^{i_n}$$
Hence by Inequality \eqref{convergence-norme},
the series \eqref{formal-series} converges normally, so $F(u_1,\cdots,u_n)\in
\boW$ and $F$ is of class $C^\infty$ on $B({\bf a},{\bf r}')$.
(See Theorem 11.12 in
\cite{chae} for the smoothness of maps defined by power series in Banach spaces.)
\cqfd

\noindent
Martin Traizet\\
Institut Denis Poisson\\
Universit\'e de Tours, 37200 Tours, France\\
\verb$martin.traizet@univ-tours.fr$

\end{document}